\documentclass[12pt,a4paper]{article}

\usepackage[dvips]{graphics}

\usepackage[dvips]{epsfig}

\newtheorem{lemma}{Lemma}

\newtheorem{cor}[lemma]{Corollary}
\newtheorem{prop}[lemma]{Proposition}

\newtheorem{rem}{Remark}
\newtheorem{conj}[lemma]{Conjecture}

\newcommand{\dimo}{\noindent \emph{Proof. }}
\newcommand{\qed}{\\ \rightline{$\Box$}\\}

\usepackage{latexsym}
\usepackage{amsfonts}
\usepackage{amssymb}
\usepackage{amsmath}
\usepackage[english]{babel}

\begin{document}

\title{A NOTE ABOUT COMPLEXITY \\ OF LENS SPACES}

\author{Maria Rita CASALI - Paola CRISTOFORI}

\maketitle

\bigskip

\begin{abstract}

Within crystallization theory, (Matveev's) complexity of a 3-manifold can be estimated by means of the combinatorial notion of {\it GM-complexity}.
In this paper, we prove that the
GM-complexity of any lens space $L(p,q)$, with $p\ge 3$, is bounded by $S(p,q)
-3$, where $S(p,q)$ denotes the sum of all partial quotients in
the expansion of $\frac q p$ as a regular continued fraction. The
above upper bound had been already established with regard to
complexity; its sharpness was conjectured by Matveev
himself
and has been recently proved for some
infinite families of lens spaces by Jaco, Rubinstein and Tillmann.
As a consequence, infinite classes of 3-manifolds turn out to exist, where
complexity and GM-complexity coincide.

Moreover, we present and briefly analyze results arising from
crystallization catalogues up to order $32$, which prompt us to
conjecture, for any lens space $L(p,q)$ with $p\ge 3$, the
following relation: $k(L(p,q)) = 5 + 2 c(L(p,q))$, where $c(M)$
denotes the complexity of a 3-manifold $M$ and $k(M) +1$
is half the minimum order of a crystallization of $M$.
\end{abstract}

\bigskip


\medskip


\bigskip
\bigskip

\section{\hskip -0.7cm . Introduction}

\medskip

The complexity $c(M)$ of a 3-manifold $M$ was originally defined
by Matveev as the minimum number of true vertices among all almost
simple spines  of $M$ (\cite{[M$_1$]}).\footnote{Recall that a
{\it spine} of $M$ is a 2-dimensional subpolyhedron $P$ of
$Int(M)$, such that $M$ (or $M$ minus an open 3-ball, if $M$ is
closed) collapses to $P$; moreover, a spine $P$ is {\it almost
simple} if the link of each point of P can be embedded in the
1-skeleton of a 3-simplex. A point whose link is the whole
1-skeleton of the 3-simplex is called a {\it true vertex} of $P$.}

The 3-sphere, the real projective space, the lens space $L(3,1)$ and the spherical bundles $\mathbb S^1\times \mathbb S^2$ and $\mathbb S^1\tilde\times \mathbb S^2$ have complexity zero. Excluding these examples, it is well-known that the complexity of a closed prime 3-manifold $M$ coincides with the number of tetrahedra in a minimal pseudocomplex triangulating $M$.
Therefore, it is possible to obtain tables of 3-manifolds for increasing values of their complexity, simply by generating all triangulations (resp. spines) with a given number of tetrahedra (resp. true vertices) and by identifying the related manifolds. 
Up-to-date censuses of closed orientable  (resp. non-orientable)  3-manifolds are available at the Web page http://www.matlas.math.csu.ru/ (resp. in \cite[Appendix]{[B$_2$]}).\footnote{See the references in \cite[section 1]{[BCCGM]} for subsequent results about the classification of 3-manifolds according to complexity, established by various authors.}

The complexity of a manifold is generally hard to compute from the theoretical point of view, leaving aside the concrete enumeration of its spines. 

As far as lens spaces are concerned, a well-known estimation for complexity 
exists: 
$c(L(p,q)) \le S(p,q) -3,$  where $S(p,q)$ denotes the sum of all partial quotients in the expansion of $\frac q p$ as a regular continued fraction 
(\cite{[M_libro]} and \cite{[JRT1]}).
The sharpness of this estimation was conjectured by Matveev himself 
(\cite[page 77]{[M_libro]});  in the last decade Jaco, Rubinstein and Tillmann succeeded to prove that it is sharp for some infinite families of lens spaces, comprehending $L(2r, 1)$ and $L(4r, 2r-1)$, $\forall r \ge 2$  (see \cite{[JRT1]} and \cite{[JRT2]}).

\medskip

In the present paper we approach the problem of computing complexity of
lens spaces by making use of  {\it crystallization theory}, i.e.
by a combinatorial tool to represent any PL-manifold
via suitable edge-coloured graphs, called {\it
crystallizations}: see \cite{[FGG]}, \cite{[BCG]} and
\cite{[BCCGM]} for a survey on the topic. Within this theory, the
invariant {\it GM-complexity} has been introduced for 3-manifolds:
see \cite{[C$_4$]} (resp. \cite {[C$_5$]} and \cite{[CC$_1$]}) for
the closed non-orientable (resp. orientable) case, and
\cite{[CC$_bordo$]} for the extension to 3-manifolds with
non-empty boundary.

The GM-complexity of a 3-manifold $M$ (denoted by $c_{GM}(M)$) turns out to be an upper bound for the complexity of $M$: via GM-complexity,
direct estimations of complexity have been obtained for interesting classes of 3-manifolds, both in the closed case (see \cite {[C$_5$]}
for 2-fold and 3-fold branched coverings of $\mathbb S^3$,  and for 3-manifolds obtained by Dehn surgery on framed links in $\mathbb S^3$)
and in the boundary case (see \cite{[CC$_bordo$]} for torus knot complements).

Conjectures about significant upper bounds for infinite families of 3-manifolds have also arisen by the analysis of existing crystallization catalogues,
which have been generated up to 30 vertices (\cite{[C$_2$]}, \cite{[CC$_2$]}, \cite{[BCrG$_1$]}).

In section \ref{catalogue_32}, we present results concerning the extension of the above catalogues to 32 vertices, with regard to the orientable case:
as a consequence, the conjectures are proved up to this order. 
In particular, we remark that
$c_{GM}(L(p,q)) = c(L(p,q)) = S(p,q) - 3$ for any lens space
involved in 
crystallization catalogues.

On the other hand, the main result of the present paper states that $S(p,q) -3$ is an upper bound
for the GM-complexity of any lens space $L(p,q)$ with $p \ge 3$ (Proposition \ref{complessità L(p,q)}(a)).

The coincidence with the known upper bound for $c(L(p,q))$ emphasizes the strong relationship between GM-complexity and complexity,
which is suggested by the analysis of the crystallization catalogues (see for example \cite[Conjecture 7]{[CCM]},
where equality between the two invariants is conjectured for all closed 3-manifolds).

In particular, lens spaces of type $L(2r, 1)$ and $L(4r, 2r-1)$, $\forall r \ge 2$, turn out to be infinite classes of 3-manifolds where complexity and GM-complexity actually coincide (see Corollary \ref{corollario main theorem}).

\medskip

Within crystallization theory, a slightly different complexity notion has been introduced and analyzed, too:
a closed PL $n$-manifold $M^n$ is said to have {\it gem-complexity} $k(M^n) = p-1$ if $2p$ is the minimum order of a crystallization representing $M^n$
(see \cite{[Li]}, \cite{[C$_2$]} and \cite{[C$_Oberwolfach$]}).

In section \ref{lens_spaces} (Proposition \ref{complessità L(p,q)} (b))
we prove that   $k(L(p,q)) \le 2 \, S(p,q) -1$ $ \forall p \ge 2$.
The sharpness of the above bound for a double infinite family of lens spaces (Corollary \ref{corollario Swartz}) and for all lens spaces up to complexity 5 suggests the equality to hold in the general case (Conjecture \ref{congettura gem-complexity lens}).

Moreover, we prove that, whenever
$c(L(p,q)) = S(p,q) -3$, then $k(L(p,q)) \le 5 + 2
c(L(p,q))$ follows  (Proposition
\ref{gem-complexity/Matveev's conjecture}). Actually, in Proposition
\ref{gem-complexity/complexity}, we list some infinite families of
lens spaces for which equality holds. Therefore, it is natural to
conjecture it for all lens spaces, i.e. \ $k(L(p,q)) = 5
+ 2 c(L(p,q))$ \ $\forall p \ge 3$  (Conjecture \ref{congettura
complexity/gem-complexity}).

\medskip
\section{\hskip -0.7cm .  Crystallizations, gem-complexity and GM-complexity}
\label{preliminari}

\medskip

Edge-coloured graphs are a representation tool for the whole class of piecewise linear (PL)
manifolds, without restrictions about dimension, connectedness, orientability or boundary properties. In the present work,
however, we will deal only with closed and connected PL-manifolds of dimension $n=3$; hence, we will briefly review basic notions and results of the theory with respect to this particular case.

A {\it 4-coloured graph (without boundary)} is a pair
$(\Gamma,\gamma)$, where $\Gamma= (V(\Gamma),E(\Gamma))$ is a
regular multigraph (i.e. it may include multiple edges, but no
loop) of degree four and $\gamma : E(\Gamma) \to
\Delta_3=\{0,1,2,3\}$ is a proper edge-coloration (i.e. it is
injective when restricted to the set of edges incident to any vertex of $\Gamma$).

\smallskip

The elements of the set $\Delta_3$ are called the {\it colours} of
$\Gamma$; thus, for every $i\in \Delta_3$, an {\it $i$-coloured
edge} is an element $e \in E(\Gamma)$ such that $\gamma(e)=i.$ For
every $i,j \in \Delta_3$ let $\Gamma_{\hat \imath}$ (resp.
$\Gamma_{ij}$) (resp. $\Gamma_{\hat \imath\hat \jmath}$) be the
subgraph obtained from $(\Gamma, \gamma)$ by deleting all the
edges of colour $i$ (resp. $c\in \Delta_3-\{i,j\}$) (resp. $c\in
\{i,j\}$). The connected components of $\Gamma_{ij}$ (resp.
$\Gamma_{\hat\imath}$) (resp. $\Gamma_{\hat\imath \hat\jmath}$)
are called {\it $\{i,j\}$-coloured} cycles (resp.
${\hat\imath}$-residues) (resp. {\it $\{\hat\imath,
\hat\jmath\}$-coloured} cycles)   of $\Gamma,$ and their number is
denoted by $g_{ij}$ (resp. $g_{\hat\imath}$) (resp. $g_{\hat\imath
\hat\jmath}$). A 4-coloured graph $(\Gamma, \gamma)$ is called
{\it contracted} iff, for each $i\in \Delta_3$, the subgraph
$\Gamma_{\hat\imath}$ is connected (i.e. iff $g_{\hat\imath}=1$
$\forall i\in \Delta_3$).

Every 4-coloured graph $(\Gamma, \gamma)$ may be thought of as the
combinatorial visualization of a $3$-dimensional labelled
pseudocomplex $K(\Gamma)$, which is constructed
according to the following instructions:
\begin{itemize}
\item{}
\  for each vertex $v\in V(\Gamma)$, take a $3$-simplex
$\sigma(v)$, with its vertices labelled $0,1,2,3$;
\item{}
\ for each $j$-coloured edge between $v$ and $w$ ($v,w\in
V(\Gamma)$), identify the 2-dimensional  faces of $\sigma(v)$ and
$\sigma(w)$ opposite to the vertex labelled $j$, so that equally
labelled vertices coincide.
\end{itemize}

In case $K(\Gamma)$ triangulates a (closed) PL 3-manifold $M$,
then $(\Gamma,\gamma)$ is called  a {\it gem} (gem = \underbar
graph  \underbar encoded \underbar manifold) {\it representing}
$M$.\footnote{The construction of $K(\Gamma)$ directly ensures
that, if $(\Gamma,\gamma)$ is a gem of $M$, then $M$ is orientable
(resp. non-orientable) iff $\Gamma$ is bipartite (resp.
non-bipartite).}

Finally, a 4-coloured graph representing a (closed) 3-manifold $M$
is  a {\it crystallization} of $M$ if it is also a contracted
graph; by construction, it is not difficult to check that this is
equivalent to requiring  that the associated pseudocomplex
$K(\Gamma)$ contains exactly one $i$-labelled vertex, for every
$i\in \Delta_3$. The representation theory of PL-manifolds by
edge-coloured graphs is often called {\it crystallization theory},
since it has been proved that every PL-manifold admits a
crystallization: see Pezzana Theorem and its subsequent
improvements (\cite{[FGG]} or \cite{[BCG]}).

\bigskip

A cellular  embedding of a coloured graph into a surface is said
to be {\it regular} if its regions are bounded by the images of
bicoloured cycles; interesting results of crystallization theory
(mainly related to an $n$-dimensional extension of Heegaard genus,
called {\it regular genus} and introduced in \cite{[G$_2$]}) rely
on the existence of this type of embeddings for graphs
representing manifolds of arbitrary dimension.

As far as the 3-dimensional case is concerned, it is well-known
that, if $(\Gamma, \gamma)$ is a crystallization of an orientable
(resp. non-orientable) manifold $M,$ then for every pair $\alpha,
\beta\in \Delta_3$, there exists a regular embedding $i_{\alpha,
\beta}\ :\ \Gamma \to F_{\alpha \beta}$, such that $F_{\alpha
\beta}$ is the closed orientable  (resp. non-orientable) surface
of genus $g_{\alpha\beta} - 1.$ The minimum genus of $F_{\alpha
\beta}$ taken over all pairs $\alpha, \beta\in \Delta_3$ is called
the {\it regular genus} of $(\Gamma, \gamma)$.

\bigskip

Let now $\mathcal D$ (resp. $\mathcal
D^{\prime}$) be an arbitrarily chosen $\{\alpha, \beta\}$-coloured
(resp. $\{\hat \alpha, \hat \beta\}$-coloured) cycle of a crystallization $(\Gamma,
\gamma);$ we denote by $\mathcal R_{\mathcal D, \mathcal
D^{\prime}}$ the set of regions of $F_{\alpha \beta}-i_{\alpha, \beta}((\Gamma_{\alpha
\beta} - \mathcal D)\cup (\Gamma_{\hat \alpha \hat \beta} -
\mathcal D^{\prime})).$

\bigskip

\par \noindent
{\bf Definition 1.} Let $M$ be a closed 3-manifold, and let
$(\Gamma, \gamma)$ be a crystallization of $M.$  With the above
notations, the {\it Gem-Matveev complexity} (or {\it GM-complexity},
for short) of $\Gamma$ is defined as the non-negative integer
$$ c_{GM}(\Gamma) = min \{\# V(\Gamma)-\# (V(\mathcal D)\cup V(\mathcal D^{\prime}) \cup V(\Xi)) \ /
\ \mathcal D \in \Gamma_{\alpha \beta}, \mathcal D^{\prime} \in
\Gamma_{\hat \alpha \hat \beta}, \Xi \in \mathcal R_{\mathcal D,
\mathcal D^{\prime}}\},$$ while the {\it Gem-Matveev complexity}
(or {\it GM-complexity}, for short) of $M$ is defined as
$$ c_{GM}(M) = min \{c_{GM}(\Gamma) \ / \ (\Gamma, \gamma) \
\text{crystallization of} \ M \}$$

\smallskip

\begin{rem}
{\em The notion of  GM-complexity can also be extended to
non-contracted gems (see \cite[Definition 4]{[C$_5$]}); however,
it has been proved that, for each closed 3-manifold $M$, the
minimum value of GM-complexity is always realized by
crystallizations (\cite[Proposition 6]{[CCM]}). For this reason,
in this paper we restrict our attention to GM-complexity of
crystallizations.}
\end{rem}

The following key result, due to \cite{[C$_4$]}, justifies the choice of
terminology:

\begin{prop}\label{relazione c/cGM}
For every closed 3-manifold $M$, GM-complexity gives an
upper bound for complexity of $M$:
$$c(M) \le c_{GM}(M)$$
\vskip -1truecm \ \qed
\end{prop}

Within crystallization theory, another quite natural notion of
complexity has been introduced (\cite{[Li]}, \cite{[C$_2$]}):

\bigskip

\par \noindent
{\bf Definition 2.} For each closed PL $n$-manifold $M^n$,
{\it gem-complexity} is defined as the non-negative integer
$$ k(M^n) = min \{ \frac {\# V(\Gamma)}2 -1 \ / \ (\Gamma, \gamma) \
\text{crystallization of} \ M^n \}$$

\medskip

As far as dimension three is concerned, existing 3-manifold censuses via crystallizations (see paragraph \ref{catalogue_32}) induce to investigate possible relations between complexity
$c(M)$ and gem-complexity $k(M)$. In particular, in \cite{[CC$_1$]}, the following conjecture is stated:

\begin{conj}
\label{congettura gem-complexity}
\ $ k(M) \le 5 +2 c(M)$ \ \ for any closed orientable
3-manifold $M.$
\end{conj}

\medskip
\section{\hskip -0.7cm .  Up-to-date results from crystallization catalogues}
\label{catalogue_32}

\medskip

The totally combinatorial nature of graphs encoding manifolds,
makes the theory particularly effective for generating catalogues
of PL-manifolds for increasing order of the representing graphs (i.e. for increasing values of $k(M^n)$). The
main tool for the algorithmic generation of tables of
crystallizations is the {\it code}, a numerical ``string" which
completely describes the combinatorial structure of a coloured
graph, up to colour-isomorphisms (see \cite{[CG$_2$]});
afterwards, suitable moves on gems (see \cite{[FG$_1$]}) are
applied, to  
yield
a classification procedure which allows to
detect crystallizations of the same manifold.

In particular, the generation and classification procedures have
been successfully de\-ve\-loped in dimension 3, both in the orientable
and non-orientable case, allowing the complete topological
identification of each involved 3-manifold: see \cite{[Li]}
and \cite{[CC$_2$]} (resp. \cite{[BCrG$_1$]}) for  censuses of
orientable (resp. non-orientable) 3-manifolds up to gem-complexity
$14$.

\medskip

In the orientable case (to which we restrict our attention in the
present paper), the obtained results 
allow to state:
 
\begin{prop}
There are exactly $110$ closed prime orientable
$3$-manifolds, up to gem-complexity $14$.
\end{prop}

\noindent

Details about the quoted catalogues are available at the Web page:
\par \noindent
http://cdm.unimo.it/home/matematica/casali.mariarita/CATALOGUES.htm.
\par
In particular, Table 1,  Table 2 and  Table 3 of the above WEB
page contain the JSJ-decomposition of each manifold up to
gem-complexity $14$ and an analysis of the distribution of these
manifolds with respect to complexity and geometry.

\bigskip

By optimizing the program code and by exploiting high-powered
computers\footnote{In particular, we made use of CINECA
facilities, such as a IBM SP Power6 system for high-performance
computing, which are available in virtue of some established
Italian Supercomputing Resource Allocation (ISCRA) projects.}, we
succeeded in extending the above crystallization catalogue to
3-manifolds with $k(M)=15$.\footnote{The same task has been
independently achieved by Tarkaev and Fominykh, too: see
\cite[page 367]{[M_libro]}.}

The classification of the involved manifolds has been performed by
means of the powerful computer program {\it ``3-Manifold
Recognizer''}.\footnote{``3-Manifold Recognizer'' has been written
by V. Tarkaev as an application of the results about recognition
of 3-manifolds obtained by S. Matveev and his research group. It
is available on the Web: http://www.matlas.math.csu.ru/}

\medskip

The analysis of the catalogue yields the following result, concerning the geometries of the involved manifolds
and the comparison between $k(M)$ and $c(M)$:

\begin{prop} \ There \! are \! exactly \! 110 \! closed \! prime
\! orientable \! 3-manifolds \! with \! gem-complexity $15$. They
are:
\begin{itemize}
\item{} 38  elliptic 3-manifolds (in particular: 20 lens spaces,
which are exactly those with complexity five; $18$ other elliptic,
of which $11$ with complexity six \footnote{They are exactly the
missing ones with that complexity in the crystallization catalogue
for  $k \le 14$.} and seven with complexity seven); \item{} four
Nil Seifert 3-manifolds (one with complexity seven, one with
complexity eight and two with complexity nine); 
\item{} ten
3-manifolds with Sol geometry (one with complexity seven \footnote{It is exactly the missing one with that complexity in
the crystallization catalogue for $k \le 14$.}, three with
complexity eight and six with complexity nine); 
\item{} two
3-manifolds with $\mathbb H^2 \times \mathbb R$ geometry (which
are exactly those with complexity eight); \item{} 34 3-manifolds
with $\widetilde{SL}_2(\mathbb R)$ geometry (six with complexity
seven, $19$ with complexity eight and nine with complexity nine);
\item{} 17 non-geometric 3-manifolds (five with complexity eight
and $12$ with complexity nine); \item{} five hyperbolic
3-manifolds  (one with complexity nine, three with complexity ten
and one with complexity $11$).\footnote{They are all Dehn fillings
of the complement of $6^3_1$ (i.e. the chain link with three
components in $\mathbb S^3$).}
\end{itemize}
 \end{prop}

\bigskip

Apart from $\mathbb S^3$, $L(2,1)$, $L(3,1)$ and the sphere bundle
$\mathbb S^2 \times \mathbb S^1$ (i.e. the ones with complexity zero), the distribution of prime 3-manifolds with
$k(M)\le 15$ with respect to complexity and geometry is summarized
in the following table (where the symbol $x/n$ means that $x$
3-manifolds appear in the catalogue among the $n$ ones having the
corresponding complexity and geometry, and bold character is used
to indicate the presence of all manifolds of the considered type).

\bigskip

\centerline{{\scriptsize{\begin{tabular}{r|c c c c c c c c c c c }{complexity} & 1  & 2 & 3& 4& 5 & 6 & 7 & 8 & 9 & 10 & 11 \\
\hline
\hline
lens & {\bf 2/2} &  {\bf 3/3} &  {\bf 6/6} &  {\bf 10/10} &  {\bf 20/20} & 0/36 & 0/72 & 0/136 & 0/272 & 0/528 & 0/1056 \\
other elliptic & - &  {\bf 1/1} &  {\bf 1/1} &  {\bf 4/4} &  {\bf 11/11} & {\bf 25/25} & 7/45 & 0/78 & 0/142 & 0/270 & 0/526\\
$\mathbb E^3$ &  - & - & - & - & - &  {\bf 6/6} & - & - & - & - & -  \\
Nil       &  - & - & - & - & - &  {\bf 7/7} & 4/10 & 3/14 & 2/15 & 0/15 & 0/15 \\
$\widetilde{SL}_2(\mathbb R)$  & - & - & - & - & - & - & 19/39 & 24/162 & 11/513 &
0/1416 & 0/3696 \\ Sol  & - & - & - & - & - & - & {\bf 5/5} & 5/9 & 6/23 & 0/39 & 0/83 \\
non-geometric & - & - & - & - & - & - &  {\bf 4/4} & 6/35 & 14/185 & 0/777 & 0/2921 \\
$\mathbb H^2 \times \mathbb R$  & - & - & - & - & - & - & - & {\bf 2/2} & - & 0/8 & 0/4 \\
hyperbolic  & - & - & - & - & - & - & - & - &  3/4 & 4/25 & 1/120 \\
\hline  \hline TOTAL &  {\bf 2/2} &  {\bf 4/4} &  {\bf 7/7} &  {\bf 14/14} &
{\bf 31/31} & 38/74 & 39/175 & 40/436 & 36/1154 & 4/3078 & 1/8421
\end{tabular}}}}

\medskip
\centerline{\footnotesize{Prime 3-manifolds with complexity $c \ne 0$ and    
gem-complexity 
$k \le 15$.}}

\bigskip

It is worthwhile to note that this new table emphasizes an
idea already suggested by the cases up to gem-complexity $14$: for
any fixed $c$, subsequent crystallization catalogues,
for increasing order, appear to cover first the  most
``complicated" types of complexity $c$ 3-manifolds and then the
simplest ones, such as lens spaces (see the above quoted Table 3).

\medskip

On the other hand 3-manifold censuses via crystallizations suggest that ``restricted" gem-complexity
implies ``restricted" complexity, and viceversa: see \cite[Section
5]{[C$_2$]}, \cite[Remark 1]{[C$_4$]} and \cite[Propositions 5 and
6]{[CC$_1$]}).
In particular, if we improve the known results in the orientable
case with the ones obtained from the analysis of the catalogue of
order $32$ bipartite crystallizations, we can state:\footnote{Note that the results of Proposition \ref{relazione c/k} strengthen Conjecture \ref{congettura gem-complexity}.}

\begin{prop} \label{relazione c/k}
Let $M$ be a closed orientable prime 3-manifold.
\begin{itemize} \item[(a)]  If $c(M) \le 6$ and $M \ne
L(p,q)$, then  $ k(M) <  5 + 2c(M);$
\item[(b)]  if $M =L(p,q)$ and  $1 \le c(M) \le 5$, then $ k(M) = 5 + 2c (M).$
\end{itemize}
\end{prop}

\vskip-1truecm
\ \qed

\medskip
\section{\hskip -0.7cm .  Results about lens spaces}
\label{lens_spaces}

\medskip

We are now able to state our main result:

\begin{prop}[Main Theorem]  \label{complessità L(p,q)}
Let $L(p,q)$ be a lens space. Then:
\begin{itemize}
\item[(a)] $c_{GM}(L(p,q)) \le S(p,q) -3$  \ \ $\forall p \ge 3;$
\item[(b)] $k(L(p,q)) \le 2 \cdot S(p,q) -1$  \ \ $\forall p \ge 2.$
\end{itemize}
Moreover, if $c(L(p,q)) \le 5,$ then:
$$k(L(p,q)) = 2 \cdot S(p,q) -1 \ \ \ \ \text{and} \ \ \ \ c_{GM}(L(p,q))= c(L(p,q))= S(p,q)-3.$$
\end{prop}

\dimo Without loss of generality, we may assume $p,q \in \mathbb Z$ coprime with $ 1 \le q \le \frac p 2$. 
Let $\mathfrak b (p,q)$ be the 2-bridge knot or link, with numbers $p$ and $q$ in Schubert's normal form, so that $L(p,q)$ is the 2-fold covering of $\mathbb S^3$ branched over $\mathfrak b (p,q)$ (see \cite[Proposition 12.3]{[BZ]}).

According to \cite[Proposition 12.13]{[BZ]}, $\mathfrak b (p,q)$ admits a presentation $\bar P$ as a 4-plat with a defining braid
$\mathfrak z= \sigma_2^{a_1} \cdot \sigma_1^{-a_2} \dots \sigma_1^{-a_{m-1}} \cdot \sigma_2^{a_m}$, with $a_i >0$ $\forall i \in \{1, \dots, m\}$
and $m$ odd, where the $a_i$'s are the quotients of the continued fraction $[a_1, \dots , a_m]= \frac q p$.
In Fig. 1 the presentation $\bar P$ is shown, for the case $p=21,\ q=8:$ in fact,  $ \frac 8 {21}=[2,1,1,1,2].$
Note that $\bar P$ admits $\sum_{i=1}^m a_i = S(p,q)$ crossings; for sake of notational simplicity, we set $S(p,q)=\bar s$.

\smallskip
\centerline{\scalebox{1}{\includegraphics{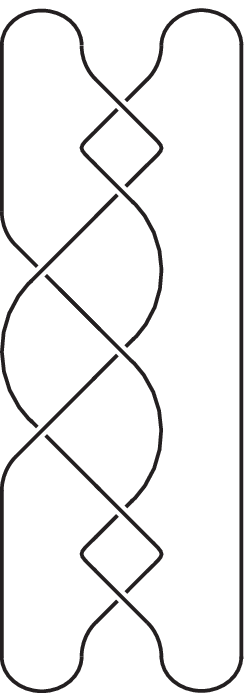}}}
\centerline{\footnotesize{Figure 1: a 4-plat presentation of $\mathfrak b (21,8)$}}

\medskip

Let now $(\bar \Gamma, \bar \gamma) = F(\bar P)$ be the
crystallization obtained by means of Ferri's construction of
crystallizations of 2-fold branched coverings (see \cite{[F]}),
 applied to $\bar P$ considered as a bridge presentation of $\mathfrak b (p,q)$, with a bridge for each crossing point.\footnote{Crystallizations
 arising from Ferri's construction are said to be {\it 2-symmetric} (see \cite{[CGr_2]}) since an involution exists on their edges,
 which exchanges colour $0$ (resp. $2$) with colour $1$ (resp. $3$): this involution can be thought of as the planar symmetry whose axis contains
 all the bridges of the given knot or link presentation.}

$(\bar \Gamma, \bar \gamma)$ may be easily ``drawn over" $\bar P$,
as it can be seen in Fig. 2 for the case $p=21,\ q=8,$ starting from the bridge presentation of the 4-plat presentation of $\mathfrak b (21,8)$ depicted in Fig. 1.

In this way,
each $\{0,1\}$-coloured cycle of $\bar \Gamma$ arises from a
bridge and has length four; hence, $(\bar \Gamma, \bar \gamma)$
turns out to have $4 \bar s$ vertices, which directly proves part
(b) of the statement.

\bigskip

\centerline{\scalebox{1}{\includegraphics{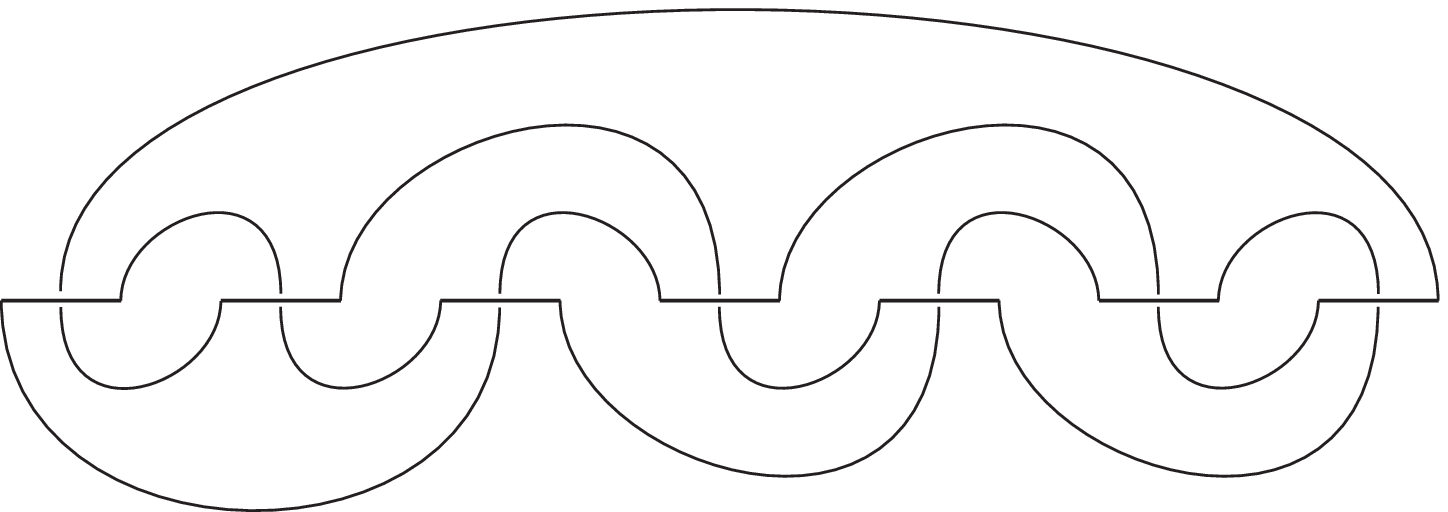}}}

\medskip

\centerline{\scalebox{1}{\includegraphics{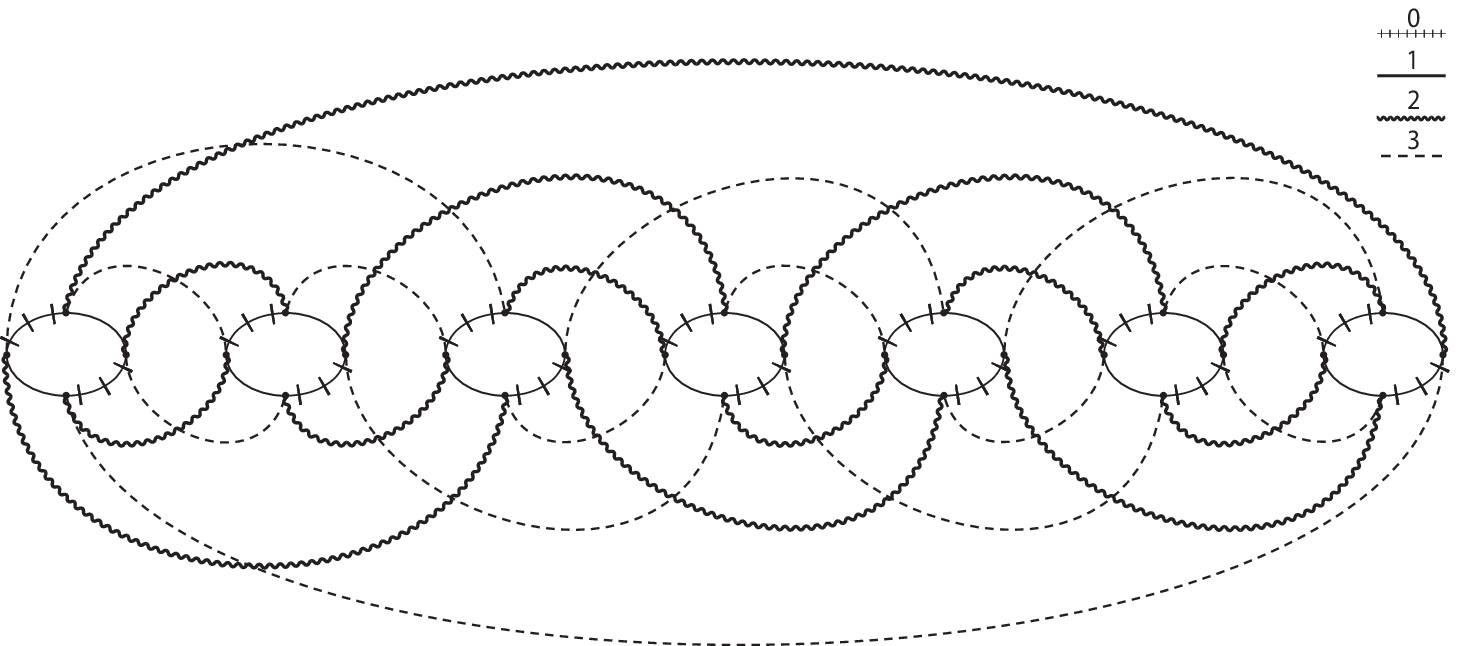}}}
\centerline{\footnotesize{Figure 2: a bridge presentation of $\mathfrak b (21,8)$ and the associated crystallization of $L(21,8)$}}
\medskip

Suppose now $p \ge 3$, i.e. $\bar s\ge 3$. Let us label the crossings of $\bar P$ according to the braid sequence
$\sigma_2^{a_1} \cdot \sigma_1^{-a_2} \dots \sigma_1^{-a_{m-1}} \cdot \sigma_2^{a_m}$ and
let $\pi$  be the ``ideal" axis containing all bridges of $\bar P$ with the orientation induced by the sequence $a_1,\ldots,a_m$.
For each crossing $c_j$ ($1\le j\le \bar s$), we label $v_{j,1}, v_{j,2}, v_{j,3}, v_{j,4}$ the vertices of its associated $\{0,1\}$-coloured cycle,
so that $v_{j,2}$ and  $v_{j,4}$ belong to $\pi$ and the sequence $v_{1,2},v_{1,4}, v_{2,2},v_{2,4}, \dots, v_{\bar s,2},v_{\bar s,4}$ is consistent with the fixed orientation of $\pi$. Moreover, $v_{j,1}, v_{j,2}$ and $v_{j,3}, v_{j,4}$ (resp. $v_{j,2}, v_{j,3}$ and $v_{j,4}, v_{j,1}$) are 0-adjacent (resp. 1-adjacent), for each $j \in 1, \dots, \bar s.$
The 2-coloured edges follow the arcs of $\bar P$, while 3-coloured edges are obtained from arcs of $\bar P$ via symmetry with respect to the axis $\pi$ (see again Fig. 2).

As a consequence we have that, for each $j,k\in\{1,\ldots, \bar s\}$, vertices $v_{j,2}$ and $v_{k,1}$ (resp. $v_{j,3}$ and $v_{k,4}$) are 3-adjacent iff $v_{j,2}$ and $v_{k,3}$ (resp. $v_{j,1}$ and $v_{k,4}$) are 2-adjacent.

In order to estimate $c_{GM}(\bar \Gamma)$ (via Definition 1), let
us now fix $\{\alpha, \beta\}=\{0,2\}$ and $\{\alpha^{\prime},
\beta^{\prime}\}=\{1,3\}$, and let us consider the
$\{0,2\}$-coloured  (resp. $\{1,3\}$-coloured) cycle $\mathcal D$
(resp. $\mathcal D^{\prime}$) corresponding to the boundary of the
``external" region $\bar R$ of $\bar P$ (resp. obtained via
symmetry from the boundary of the ``external" region $\bar R$ of
$\bar P$), and the region $\Xi \in \mathcal R_{\mathcal D,
\mathcal D^{\prime}}$ induced by the inner region $\bar
R^{\prime}$ of $\bar P$ containing the fourth string of the
4-braid associated to the 4-plat $\bar P$.

Let us define the following sets of integers:

$$I_1 = \{j\in\{2,\ldots, \bar s-1\}\ |\ c_j\text{ corresponds to a factor } \sigma_1 \text{ in } \mathfrak z\}\cup\{\bar s\}$$

$$I_2 = \{j\in\{1,\ldots,\bar s-1\}\ |\ c_j\text{ corresponds to a factor } \sigma_2  \text{ in } \mathfrak z\}.$$

Since the boundary of $\bar R$ contains the crossing points $c_j$
with $j\in I_1\cup\{1\}$, it is not difficult to check that
$\mathcal D \cup \mathcal D^{\prime}$ consists of the vertices
$\{v_{j,i} \ / \  j \in I_1,\ i = 1,3,4\} \cup \{v_{1,1}, v_{1,2},
v_{1,3}\}.$

On the other hand, the structure of $\bar P$ and the symmetry properties of $\bar\Gamma$, which arise from Ferri's construction,
allow to check that the region $\Xi$ can be obtained by the union of the four regions bounded by the following cycles of $\bar\Gamma$:
\begin{itemize}
\item [-] the $\{1,2\}$-coloured  cycle corresponding to $\bar R^{\prime}$; its vertex-set is
$\{v_{j,i}\ |\ j\in I_2,\ i=1,4\}\cup\{v_{\bar s,1},v_{\bar s,4}\}$,
    \item [-] the $\{0,3\}$-coloured cycle obtained via symmetry from $\bar R^{\prime}$; its vertex-set is
    $\{v_{j,i}\ |\ j\in I_2,\ i=3,4\}\cup\{v_{\bar s,3},v_{\bar s,4}\}$,
        \item [-] the $\{0,1\}$-coloured cycle corresponding to the last crossing $c_{\bar s}$,
        \item [-] the $\{2,3\}$-coloured cycle containing the edge corresponding to the fourth string of the braid $\mathfrak z$
        (whose vertices are exactly $\{v_{1,1}, v_{1,3}, v_{2,2}, v_{\bar s,4}\}$, since $q < \frac p 2$ implies $a_1 \ge 2$).
\end{itemize}

Hence,  $V(\mathcal D)\cup V(\mathcal D^{\prime}) \cup V(\Xi)$ contains all vertices $v_{j,i}$ with $1\le j\le \bar s,\ \ i\in\{1,3,4\}$
and the vertices $v_{j,2}$ with $j\in\{1,2,\bar s\}$.
As a consequence,
$$   V(\bar \Gamma)- (V(\mathcal D)\cup V(\mathcal D^{\prime}) \cup
V(\Xi))  =  \{ v_{j,2} \ / \  3 \le j \le \bar s-1\},$$ which
implies $ c_{GM}(\bar \Gamma)  \le   \bar s-3.$ So, part (a) of
the statement directly follows.

Finally, the last part of the statement can be checked by means of the existing crystallization catalogues, up to order 32 (presented in the previous paragraph). \qed

\medskip

Proposition \ref{complessità L(p,q)} (a), together with the results
of \cite{[JRT1]} and \cite{[JRT2]}, directly proves the existence
of infinite classes of 3-manifolds where complexity and
GM-complexity coincide:

\begin{cor}\label{corollario main theorem}\ \
\begin{itemize}
\item $c_{GM}(L(2r, 1)) = c(L(2r, 1)) = 2r-3$, for each $r \ge 2;$
\item $c_{GM}(L(4r, 2r-1)) = c(L(4r, 2r-1))=r$, for each $r \ge 2;$
\item  $c_{GM}(L((r+2)(t+1)+1, t+1)) = c(L((r+2)(t+1)+1, t+1)) = r+t$, for each $t > r > 1$, $r$ odd and $t$ even;
\item  $c_{GM}(L((r+1)(t+2)+1, t+2)) = c(L((r+1)(t+2)+1, t+2)) = r+t$, for each  $t > r > 1$, $r$ even and $t$ odd.
\end{itemize}
\end{cor}
\vskip-1truecm
\ \qed
  
\vfill \eject

With regard to gem-complexity, Proposition \ref{complessità L(p,q)}(b) suggests the following conjecture:

\begin{conj}  \label{congettura gem-complexity lens}
$k(L(p,q)) = 2 \cdot S(p,q) -1$  \ \ $\forall p \ge 2.$
\end{conj}

Actually, a recent result by Swartz (\cite{[Sw]}) \footnote{The preprint \cite{[Sw]} was posted in ArXiv (http://arxiv.org/abs/1310.1991v1)
a few days after we posted the first version of the present paper (http://arxiv.org/abs/1309.5728v1). It makes use of our Proposition \ref{complessità L(p,q)}(b) in order to prove
that $4(q+r)$ is exactly the minimum order of a crystallization of $L(qr+1,q)$, $\forall r,q\ge 1$ odd.}
proves that, for each $q,r\ge 1$ odd, $k(L(qr+1,q)) = 2(r+q)-1$;
since $S(qr+1,q)=q+r$ trivially holds $\forall r,q$, we can state: 

\begin{cor}\label{corollario Swartz}\ \
Conjecture \ref{congettura gem-complexity lens} is true for any lens space $L(qr+1, q)$, with $r,q\ge 1$ odd.
\end{cor}
\vskip-1truecm
\ \qed

\begin{rem}
{\em Note that our Main Theorem yields $k(L(p,1)) \le
2p -1.$ Actually, in this case, the described crystallization $(\bar \Gamma,
\bar \gamma)$  - which, by \cite{[Sw]}, realizes gem-complexity of $L(p,1)$ in the case of $p$ even - is
the well-known standard crystallization of $L(p,1)$, having $4 p$
vertices and regular genus one. On the contrary, if $q \ne 1$ is assumed, the standard
crystallization of $L(p,q)$ has regular genus one and $4 p$
vertices, too, but it does not realize minimum order, since
it contains {\it clusters} or, more generally, structures that can be eliminated (see  
\cite[section 6]{[CC$_2$]} or \cite[section 5]{[FG$_1$]} 
for details).}
\end{rem}

\begin{rem}
{\em Proposition \ref{complessità L(p,q)}(b) yields a general result, which admits \cite[Theorem 2.6(v)]{[BD]}  and  \cite[Theorem 2.6(vi)]{[BD]}
as particular cases.
Note that, in the cited paper, the authors' approach via
fundamental groups     allows to rediscover examples of minimal
crystallizations already presented in \cite{[Li]} (see
http://cdm.unimo.it/home/matematica/casali.mariarita/Table1.pdf).}
\end{rem}

Finally, we can summarize the relation between complexity and gem-complexity of lens spaces as follows:

\begin{prop}  \label{gem-complexity/complexity}
\ \
\begin{itemize}
\item[(a)]
Let $L(p,q)$ be a lens space, with $p \ge 3$. Then
$$k(L(p,q)) = 5 + 2 \cdot c(L(p,q))$$
whenever:
\begin{itemize}
\item[($a_1$)] $c(L(p,q)) \le 5;$
\item[($a_2$)] $p$ is even and
$q=1;$
\item[($a_3$)] $p=(r+2)(t+1)+1$ and $q=t+1$, for $t > r >
1$, $r$ odd and $t$ even;
\item[($a_4$)] $p=(r+1)(t+2)+1$ and
$q=t+2$, for $t > r > 1$, $r$ even and $t$ odd.
\end{itemize}
\item[(b)]
For each  $r,q \ge 1$ odd, then \
$$k(L(qr+1,q)) \ge 5 + 2 \cdot c(L(qr+1,q));$$
\item[(c)] for each  $r \ge 2$, then \
$$k(L(4r,2r-1)) \le 5 + 2 \cdot c(L(4r,2r-1)).$$
\end{itemize}
\end{prop}

\dimo Part ($a_1$), already stated in Proposition \ref{complessità
L(p,q)}(a), directly follows from the analysis of crystallization
catalogues up to gem-complexity $15.$

In order to prove part ($a_2$) (resp. ($a_3$)) (resp. ($a_4$)),
the equality $k(L(qr+1,q))= 2 (r+q)-1$ (which holds $\forall r,q
\ge 1$ odd: see Corollary \ref{corollario Swartz}) has to be
applied to the involved infinite class of lens spaces, together
with \cite[Theorem 1]{[JRT1]} (resp. \cite[Theorem 2]{[JRT1]})
(resp. \cite[Theorem 3]{[JRT1]}).

Part (b) follows from the equality $k(L(qr+1,q)= 2 (r+q)-1,$
$\forall r,q \ge 1$ odd, too, by making use of the 
well-known
general estimation $c(L(p,q)) \le S(p,q) -3.$

Part (c) is a direct consequence of Proposition \ref{complessità
L(p,q)}(b) and of \cite[Corollary 3]{[JRT2]}. \qed

Moreover, we can state

\begin{prop}  \label{gem-complexity/Matveev's conjecture}
Assuming Matveev's conjecture $c(L(p,q)) = S(p,q) -3$ to be true,
 then $$k(L(p,q)) \le 5 + 2 \cdot c(L(p,q)) \ \ \ \forall p \ge 2.$$

In particular, $k(L(p,q)) \le 5 + 2 \cdot c(L(p,q))$ for
each lens space $L(p,q)$ with $ 6 \le c(L(p,q)) \le 12$ and equality holds if $p=qr+1$, with $q,r$ odd.
\end{prop}

\dimo \ If \, $ c(L(p,q)) = S(p,q) -3 \, $  is assumed \, $ \forall p \ge 3, \, $ \ then \ $ k(L(p,q)) \le 2 \cdot S(p,q) -1 = \hfill$ \\  $5 + 2
\cdot c(L(p,q))$ trivially follows from part (b) of our Main Theorem.

In case $6\le c(L(p,q))\le 12$, condition $c(L(p,q)) = S(p,q) -3$ is always satisfied (see http://www.matlas.math.csu.ru/) and, therefore, $k(L(p,q)) \le 5 + 2 \cdot c(L(p,q))$. On the other hand, if, further, $p=qr+1$ ($q,r$ odd)
then part (b) of Proposition \ref{gem-complexity/complexity} applies to ensure equality.\qed

The above Propositions \ref{gem-complexity/complexity} and \ref{gem-complexity/Matveev's conjecture} naturally
suggest the following:

\begin{conj}  \label{congettura complexity/gem-complexity}
$k(L(p,q)) = 5 + 2 \cdot c(L(p,q))$  \ \ $\forall p \ge 3.$
\end{conj}

\begin{rem}
{\em A slightly different notion of GM-complexity for closed
3-manifolds was originally defined by taking into account only the
values of $c_{GM}(\Gamma),$ when $(\Gamma, \gamma)$ is a
\emph{minimal} crystallization with respect to the number of
vertices: see \cite[Definition 3]{[C$_4$]}.  As regards lens spaces, existing crystallization catalogues
prove that this restricted notion actually coincides with the one
including \emph{all} crystallizations. Moreover, the arguments
of the present paper induce {\em to conjecture} - according to the
above Conjectures \ref{congettura gem-complexity lens} and
\ref{congettura complexity/gem-complexity} - {\em that complexity
of any lens space is always realized by GM-complexity of a minimal
crystallization} (despite what happens in the general case: see \cite[Proposition 7 and Remark 4]{[CC$_1$]}).}
\end{rem}

\smallskip

{\it Acknowledgements.} We wish to thank E.~Fominykh for the
helpful discussions which originated this paper. We are also
indebted to him, to S.~Matveev and V.~Tarkaev for the use of
``3-Manifold Recognizer''.

This work was performed under the auspices of the
G.N.S.A.G.A. of I.N.d.A.M. (Italy) and financially supported by
M.I.U.R. of Italy, University of Modena and Reggio Emilia, funds for selected research topics.

\small{}

\end{document}